\documentclass[letterpaper]{amsart}

\usepackage{amssymb}
\usepackage{amsthm}
\usepackage{amsmath}
\usepackage[all]{xy} \SelectTips{eu}{}
\usepackage[hidelinks]{hyperref}

\newcommand{\numberseries}{\mdseries}   

\newlength{\thmtopspace}                
\newlength{\thmbotspace}                
\newlength{\thmheadspace}               
\newlength{\thmindent}                  

\newtheoremstyle{fixed bf head,slanted body}
                {\thmtopspace}{\thmbotspace}{\slshape}
                {\thmindent}{\bfseries}{.}{\thmheadspace}
                {{\numberseries \thmnumber{(#2) }}\thmname{#1}\thmnote{ (#3)}}

\newtheoremstyle{fixed bf head,upright body}
                {\thmtopspace}{\thmbotspace}{\upshape}
                {\thmindent}{\bfseries}{.}{\thmheadspace}
                {{\numberseries \thmnumber{(#2) }}\thmname{#1}\thmnote{ (#3)}}

\newtheoremstyle{bfupright head,upright body}
                {\thmtopspace}{\thmbotspace}
                {\upshape}{\thmindent}{\bfseries}{.}{\thmheadspace}
                {{\numberseries \thmnumber{(#2) }}\thmnote{#3}}

\newtheoremstyle{numbered paragraph}
                {\thmtopspace}{\thmbotspace}{\upshape}
                {\thmindent}{\upshape}{}{0pt}
                {{\numberseries \thmnumber{(#2) }}}

\theoremstyle{fixed bf head,slanted body}
\newtheorem{thm}{Theorem}[section]          \newtheorem*{thm*}{Theorem}
\newtheorem{prp}[thm]{Proposition}      \newtheorem*{prp*}{Proposition}
\newtheorem{cor}[thm]{Corollary}        \newtheorem*{cor*}{Corollary}
\newtheorem{lem}[thm]{Lemma}            \newtheorem*{lem*}{Lemma}
\newtheorem{fct}[thm]{Fact}

\theoremstyle{fixed bf head,upright body}
\newtheorem{dfn}[thm]{Definition}       \newtheorem*{dfn*}{Definition}
\newtheorem{rmk}[thm]{Remark}           \newtheorem*{rmk*}{Remark}
\newtheorem{exa}[thm]{Example}          \newtheorem*{exa*}{Example}

\theoremstyle{bfupright head,upright body}
\newtheorem{bfhpg}[thm]{}               \newtheorem*{bfhpg*}{}

\theoremstyle{numbered paragraph}
\newtheorem{ipg}[thm]{}

\newlength{\thmlistleft}        
\newlength{\thmlistright}       
\newlength{\thmlistpartopsep}   
\newlength{\thmlisttopsep}      
\newlength{\thmlistparsep}      
\newlength{\thmlistitemsep}     

\newcounter{eqc} 
\newenvironment{eqc}{\begin{list}{\upshape (\textit{\roman{eqc}})}%
    {\usecounter{eqc}%
      \setlength{\leftmargin}{\thmlistleft}%
      \setlength{\labelwidth}{\thmlistleft}%
      \setlength{\rightmargin}{\thmlistright}%
      \setlength{\partopsep}{\thmlistpartopsep}%
      \setlength{\topsep}{\thmlisttopsep}%
      \setlength{\parsep}{\thmlistparsep}%
      \setlength{\itemsep}{\thmlistitemsep}}}%
  {\end{list}}%

\newcounter{prt}
\newenvironment{prt}{\begin{list}{\upshape (\alph{prt})}%
    {\usecounter{prt}%
      \setlength{\leftmargin}{\thmlistleft}%
      \setlength{\labelwidth}{\thmlistleft}%
      \setlength{\rightmargin}{\thmlistright}%
      \setlength{\partopsep}{\thmlistpartopsep}%
      \setlength{\topsep}{\thmlisttopsep}%
      \setlength{\parsep}{\thmlistparsep}%
      \setlength{\itemsep}{\thmlistitemsep}}}%
  {\end{list}}%

\newenvironment{prf*}[1][Proof]{%
  \begin{proof}[\bf #1]
    \setcounter{equation}{0}
    }
  {\end{proof}
}

\newcommand{\pgref}[1]{(\ref{#1})}
\newcommand{\thmref}[2][Theorem~]{#1\pgref{thm:#2}}
\newcommand{\corref}[2][Corollary~]{#1\pgref{cor:#2}}
\newcommand{\prpref}[2][Proposition~]{#1\pgref{prp:#2}}
\newcommand{\lemref}[2][Lemma~]{#1\pgref{lem:#2}}
\newcommand{\dfnref}[2][Definition~]{#1\pgref{dfn:#2}}
\newcommand{\exaref}[2][Example~]{#1\pgref{exa:#2}}
\newcommand{\rmkref}[2][Remark~]{#1\pgref{rmk:#2}}
\newcommand{\secref}[2][Section~]{#1\ref{sec:#2}}

\newcommand{\thmcite}[2][?]{\cite[thm.~#1]{#2}}
\newcommand{\corcite}[2][?]{\cite[cor.~#1]{#2}}
\newcommand{\prpcite}[2][?]{\cite[prop.~#1]{#2}}
\newcommand{\lemcite}[2][?]{\cite[lem.~#1]{#2}}
\newcommand{\seccite}[2][?]{\cite[sec.~#1]{#2}}
\newcommand{\defcite}[2][?]{\cite[def.~#1]{#2}}
\newcommand{\exacite}[2][?]{\cite[exa.~#1]{#2}}
\newcommand{\rmkcite}[2][?]{\cite[rmk.~#1]{#2}}

\numberwithin{equation}{thm}

\def\soft#1{\leavevmode\setbox0=\hbox{h}\dimen7=\ht0\advance \dimen7
  by-1ex\relax\if t#1\relax\rlap{\raise.6\dimen7
  \hbox{\kern.3ex\char'47}}#1\relax\else\if T#1\relax
  \rlap{\raise.5\dimen7\hbox{\kern1.3ex\char'47}}#1\relax \else\if
  d#1\relax\rlap{\raise.5\dimen7\hbox{\kern.9ex \char'47}}#1\relax\else\if
  D#1\relax\rlap{\raise.5\dimen7 \hbox{\kern1.4ex\char'47}}#1\relax\else\if
  l#1\relax \rlap{\raise.5\dimen7\hbox{\kern.4ex\char'47}}#1\relax \else\if
  L#1\relax\rlap{\raise.5\dimen7\hbox{\kern.7ex
  \char'47}}#1\relax\else\message{accent \string\soft \space #1 not
  defined!}#1\relax\fi\fi\fi\fi\fi\fi}
\def\cC{\mathcal C}
\newcommand{\gra}{\alpha}
\newcommand{\grb}{\beta}

\def\urltilda{\kern -.15em\lower .7ex\hbox{\~{}}\kern .04em} 
\newcommand{\setof}[3][\mspace{1mu}]{\{#1#2 \mid #3#1\}}
\newcommand{\kk}{\Bbbk}
\newcommand{\ZZ}{\mathbb{Z}}
\newcommand{\qtext}[1]{\quad\text{#1}\quad}
\newcommand{\qand}{\qtext{and}}
\newcommand{\deq}{\:=\:}

\newcommand{\dis}{\:\is\:}
\DeclareMathOperator*{\colim}{colim}
\newcommand{\f}{\varphi}

\newcommand{\s}{\sigma}
\newcommand{\is}{\cong}
\newcommand{\qis}{\simeq}
\renewcommand{\le}{\leqslant}
\renewcommand{\ge}{\geqslant}

\newcommand{\lra}{\longrightarrow}
\newcommand{\xla}[2][]{\xleftarrow[#1]{\;#2\;}}
\newcommand{\xra}[2][]{\xrightarrow[#1]{\;#2\;}}
\newcommand{\qla}{\xla{\;\qis\;}}
\newcommand{\qra}{\xra{\;\qis\;}}
\newcommand{\ira}{\xra{\;\is\;}}
\newcommand{\QQ}{\mathbb{Q}}
\newcommand{\QZ}{\QQ/\ZZ}
\newcommand{\Rop}{R^\circ}

\newcommand{\mapdef}[4][\rightarrow]{\nobreak{#2\colon #3 #1 #4}}
\newcommand{\dmapdef}[4][\lra]{\nobreak{#2\colon #3\:#1\:#4}}
\newcommand{\Ker}[1]{\nobreak{\operatorname{Ker}#1}}

\newcommand{\Image}[1]{\nobreak{\operatorname{Image}#1}}
\newcommand{\Cone}[1]{\nobreak{\operatorname{Cone}#1}}

\newcommand{\dif}[2][]{{\partial}^{#2}_{#1}}
\newcommand{\Bo}[2][]{\operatorname{B}_{#1}(#2)}
\newcommand{\Cy}[2][]{\operatorname{Z}_{#1}(#2)}
\newcommand{\Co}[2][]{\operatorname{C}_{#1}(#2)}
\renewcommand{\H}[2][]{\operatorname{H}_{#1}(#2)}
\newcommand{\Hom}[3][R]{\operatorname{Hom}_{#1}(#2,#3)}
\newcommand{\cathom}[3]{\operatorname{hom}_{#1}(#2,#3)}
\newcommand{\tp}[3][R]{\nobreak{#2\otimes_{#1}#3}}
\newcommand{\Zphat}{\ZZ_{(p)}^{\wedge}}

\newcommand{\fd}[2][R]{\operatorname{flat\,dim}_{#1}#2}
\newcommand{\pd}[2][R]{\operatorname{proj.dim}_{#1}#2}

\title{Pure-minimal chain complexes}

\author{Lars Winther Christensen} 

\address{L.W.C. \ Texas Tech University, Lubbock, TX 79409}

\email{lars.w.christensen@ttu.edu} 

\urladdr{http://www.math.ttu.edu/\urltilda lchriste}

\author{Peder Thompson} 

\address{P.T. \ Texas Tech University, Lubbock, TX 79409}

\curraddr{Norwegian University of Science and Technology}

\email{peder.thompson@ntnu.no}

\urladdr{https://folk.ntnu.no/pedertho/}

\thanks{L.W.C.\ was partly supported by Simons Foundation
  collaboration grant 428308.}

\date{\today} 
\date{3 October 2018} 



\begin{document}

\ \vspace{2\baselineskip}

\begin{abstract}
  We introduce a notion of {\em pure-minimality} for chain complexes
  of modules and show that it coincides with (homotopic) minimality in
  standard settings, while being a more useful notion for complexes of
  flat modules. As applications, we characterize von Neumann regular
  rings and left perfect rings.
\end{abstract}

\maketitle

\section*{Introduction}

\noindent
Given a chain complex, it is natural to ask whether it has a
``smallest'' subcomplex of the same homotopy type, as such a
subcomplex would carry all pertinent information of the ambient
complex without homotopic redundancy. Initiated by Eilenberg and
Zilber \cite{SElJAZ50} in the context of simplicial complexes, this
perspective has come to play a significant role in the homological
study of rings and modules.

Let $R$ be an associative unital ring. A chain complex $M$ of
$R$-modules is called {\em minimal} if every homotopy equivalence
$M\to M$ is an isomorphism; see Avramov and Martsinkovsky
\cite{LLAAMr02} and Roig~\cite{ARg93}.  Many homological invariants of
modules, such as their injective and projective dimension, can be read
off from minimal resolutions, provided that they exist. Minimal
injective and minimal flat resolutions exist for every $R$-module;
indeed, any resolution constructed from injective envelopes or from
flat covers is a minimal chain complex; see e.g.\ Thompson
\cite{PTh}. Over a perfect ring every module has a minimal projective
resolution, and over a semi-perfect ring every finitely generated
module has a minimal projective resolution.

While minimal flat resolutions exist, they do not quite behave as one
might hope. For example, let $p$ be a prime and consider the local
ring $\ZZ_{(p)}$ with $p\ZZ_{(p)}$-adic completion $\Zphat$. It is
elementary to verify that \mbox{$F=0\to \ZZ_{(p)} \to \Zphat \to 0$}
is a minimal chain complex of $\ZZ_{(p)}$-modules; see
\exaref{minimal_not_puremin}. Evidently, it is a flat resolution of
the module $\Zphat/\ZZ_{(p)}$; however, $\Zphat/\ZZ_{(p)}$ is a flat
$\ZZ_{(p)}$-module and, as such, a minimal flat resolution of
itself. This non-uniqueness of minimal flat resolutions is, perhaps,
unsurprising as flat resolutions do not come with comparison maps the
way projective and injective resolutions do. The difference between
the flat resolutions $F$ and \smash{$\Zphat/\ZZ_{(p)}$} is the pure
subcomplex
\mbox{$P = 0\to \ZZ_{(p)} \xrightarrow{=} \ZZ_{(p)} \to 0$}; indeed,
there is a pure exact sequence
$\smash{0\to P\to F\to \Zphat/\ZZ_{(p)} \to 0}$ of chain
complexes. This points to a notion of minimality that forbids the
existence such subcomplexes.

In this paper, a chain complex is called {\em pure-minimal} if the
zero complex is the only pure-acyclic pure subcomplex. This definition
is inspired by the example recounted above and by the fact that in a
minimal chain complex the zero complex is the only contractible direct
summand.  In settings where minimality is well-understood, such as for
chain complexes of projective modules over a perfect ring, we show
that pure-minimality coincides with minimality
(\thmref{pureminiffweakmin}).

Our central construction (\thmref{sub_quotient_complexes}) shows that
given a chain complex $M$ of $R$-modules, there exists a pure-minimal
chain complex that is isomorphic to $M$ in the derived category over
$R$. Moreover, in settings where minimal chain complexes are known to
exist, the construction recognizes them (\corref{decomposition}).

As an application of our construction, we show (\thmref{pm_fd}) that
for every chain complex $M$ of $R$-modules there exists a pure-minimal
semi-flat\footnote{\,In the literature, e.g.\ in \cite{LLAHBF91} by
  Avramov and Foxby, such complexes are also called dg-flat.}  complex
$F$ that is isomorphic to $M$ in the derived category over $R$, and
that the flat dimension of $M$ can be read off from $F$. In fact,
pure-minimality is also an appropriate notion of minimality for
degreewise finitely generated semi-projective complexes over a
noetherian ring (\thmref{pm_pd}). As further applications, we
characterize von Neumann regular rings (\corref{vnr}) and left perfect
rings (\thmref{perfect_pm}).

The paper is organized as follows. In \secref{pac} we study
pure-acyclic chain complexes and give a characterization of von
Neumann regular rings in terms of pure homological algebra
(\thmref{vnr}); in \secref{pqi} we continue with a few technical
results on pure quasi-isomorphisms. In \secref{flavors_of_min} we
define pure-minimality and compare it with other notions of
minimality. We focus separately on minimality of acyclic chain
complexes in \secref{min_ac}. The main results advertised above are
found in \secref{pm_replacements}. In the appendix we establish
sufficient conditions for acyclicity of chain complexes.
\begin{equation*}
  * \ \ * \ \ *
\end{equation*}
\noindent
Throughout, $R$ is an associative algebra over a commutative unital
ring $\kk$ which, if no other choice is more appealing, can be
$\ZZ$. The term \emph{$R$-module} refers to a left $R$-module, while a
right $R$-module is considered to be a (left) module over the opposite
ring $\Rop$.

A chain complex of $R$-modules is for short called an
\emph{$R$-complex}. The category of $R$-complexes is denoted
$\cC(R)$. For a homologically indexed $R$-complex $M$, write
$\partial^M$ for the differential and define the subcomplexes $\Cy{M}$
and $\Bo{M}$ with zero differentials by specifying their modules:
$\Cy[i]{M}=\Ker(\partial_i^M)$ and
$\Bo[i]{M} = \Image(\partial_{i+1}^M)$. Further, set $\Co{M}=M/\Bo{M}$
and $\H{M} = \Cy{M}/\Bo{M}$. A complex $M$ is said to be {\em acyclic}
if the sequence $0 \to \Cy[i]{M} \to M_i \to \Cy[i-1]{M} \to 0$ is
exact for every $i\in \ZZ$; that is, $\H{M}$ is the zero complex. A
complex $M$ is called \emph{contractible} if the identity $1^M$ is
null-homotopic; a contractible complex is acyclic.

Homology is a functor on $\cC(R)$. A morphism $\mapdef{\gra}{M}{N}$ in
$\cC(R)$ is called a \emph{quasi-isomorphism} if $\H{\gra}$ is an
isomorphism. Prominent quasi-isomorphisms are homotopy equivalences;
they are morphisms that have an inverse up to homotopy.

For $R$-complexes $L$ and $M$, the total Hom complex is denoted
$\Hom{L}{M}$. For an $\Rop$-complex $N$ and an $R$-complex $M$, the
total tensor product complex is written $\tp{N}{M}$.

\section{Pure-acyclic complexes}
\label{sec:pac}

\noindent
In this first section we recall fundamentals of pure homological
algebra, with a focus on pure-acyclicity, and give a characterization
of von Neumann regular rings.

\begin{ipg}{\bf Resolutions of complexes.}
  An $R$-complex $F$ is called \emph{semi-flat} if it consists of flat
  $R$-modules and the functor $-\otimes_RF$ preserves acyclicity. A
  complex $F$ of flat modules with $F_i=0$ for $i\ll 0$ is semi-flat,
  this follows for example from \corref{ptwize-C}. Similarly, a
  complex $I$ (a complex $P$) is called \emph{semi-injective}
  (\emph{semi-projective}) if it consists of injective modules
  (projective modules) and the functor $\Hom{-}{I}$ preserves
  acyclicity ($\Hom{P}{-}$ preserves acyclicity). Every
  semi-projective complex is semi-flat. A complex $P$ of projective
  modules with $P_i=0$ for $i\ll 0$ is semi-projective, and a complex
  $I$ of injective modules with $I_i=0$ for $i \gg 0$ is
  semi-injective, this follows from \prpref[Propositions~]{ptwize-C}
  and \prpref[]{ptwize-Z}.

  Every $R$-complex $M$ has a semi-projective resolution and a
  semi-injective resolution; that is, there are quasi-isomorphisms
  \begin{equation*}
    P \qra M \qra I
  \end{equation*}
  with $P$ semi-projective and $I$ semi-injective; see
  \cite{LLAHBF91}. For a module, a classical projective (injective)
  resolution is a semi-projective (-injective) resolution.
\end{ipg}

\begin{bfhpg}[Purity in the category of modules]
  \label{pure-m}
  An exact sequence of $R$-modules $0 \to L \to M \to N \to 0$ is
  called \emph{pure} if the induced sequence of $\kk$-modules
  \begin{equation*}
    0 \lra \Hom{A}{L} \lra \Hom{A}{M} \lra \Hom{A}{N} \lra 0
  \end{equation*}
  is exact for every finitely presented $R$-module $A$. Equivalently,
  the sequence of $\kk$-modules
  $0 \to \tp{B}{L} \to \tp{B}{M} \to \tp{B}{N} \to 0$ is exact for
  every $\Rop$-module~$B$.

  An $R$-module $F$ is flat if and only if every exact sequence
  $0 \to L \to M \to F \to 0$ is pure. For a flat $R$-module $F$, an
  exact sequence $0 \to L \to F \to N \to 0$ is pure if and only if
  $L$ and $N$ are flat. See for example Lam \seccite[4J]{Lam2} for
  details.

  On the other hand, an $R$-module $E$ is fp-injective if and only if
  every exact sequence $0 \to E \to M \to N \to 0$ is pure. For an
  fp-injective $R$-module $E$, an exact sequence
  $0 \to L \to E \to N \to 0$ is pure if and only if $L$ is
  fp-injective.
\end{bfhpg}

\begin{rmk}
  \label{rmk:vnr}
  In view of \pgref{pure-m} the following conditions are equivalent.
  \begin{eqc}
  \item Every $R$-module is flat.
  \item Every short exact sequence of $R$-modules is pure.
  \item Every $R$-module is fp-injective.
  \end{eqc}
  The rings that satisfy these conditions are precisely the von
  Neumann regular rings---also called absolutely flat rings.
\end{rmk}

\begin{dfn}
  An exact sequence of $R$-complexes $0 \to L \to M \to N \to 0$ is
  called \emph{degreewise pure} (\emph{degreewise split}) if the
  sequence $0 \to L_i \to M_i \to N_i \to 0$ of $R$-modules is pure
  (split) for every \mbox{$i\in\ZZ$}.  A subcomplex $L\subseteq M$,
  and a quotient complex $M/L$, are called degreewise pure (degreewise
  split) if the canonical exact sequence $0 \to L \to M \to M/L \to 0$
  is degreewise~pure (degreewise split).
\end{dfn}

\begin{lem}
  \label{lem:pmsm}
  Let $M$ be an $R$-complex. Under each of the conditions {\rm
    (a)--(d)} below, every exact sequence of $R$-complexes
  $0 \lra L \lra M \lra N \lra 0$ that is degreewise pure is
  degreewise split.
  \begin{prt}
  \item $R$ is left noetherian; $M$ is a complex of injective
    $R$-modules.
  \item $R$ is left perfect; $M$ is a complex of projective
    $R$-modules.
  \item $R$ is semi-perfect; $M$ is a complex of finitely generated
    projective $R$-modules.
  \item $R$ is left noetherian; $M$ is a complex of finitely generated
    projective~\mbox{$R$-modules.}
  \end{prt}
\end{lem}

\begin{prf*}
  Consider a degreewise pure exact sequence of $R$-complexes
  \begin{equation*}
    \tag{$\ast$}
    0 \lra L \lra M \lra N \lra 0\:.
  \end{equation*}
  (a): It follows from \pgref{pure-m} that $L$ is a complex of
  fp-injective modules. As $R$ is left noetherian, fp-injective
  $R$-modules are injective. Thus, $(\ast)$ is degreewise an exact
  sequence of injective modules, in particular it is degreewise split
  exact.

  (b): It follows from \pgref{pure-m} that $N$ is a complex of flat
  modules. As $R$ is left perfect, flat $R$-modules are
  projective. Thus, $(\ast)$ is degreewise an exact sequence of
  projective modules, in particular it is degreewise split exact.

  (c) \& (d): The complex $N$ is degreewise finitely generated and as
  in (b) a complex of flat modules. Over a semi-perfect or left
  noetherian ring, finitely generated flat modules are projective, so
  as in (b) the sequence $(\ast)$ is degreewise split exact.
\end{prf*}

\begin{dfn}
  An $R$-complex $M$ is called \emph{pure-acyclic} if it is acyclic
  and the exact sequence
  $0 \to \Cy[i]{M} \to M_i \to \Cy[i-1]{M} \to 0$ is pure for every
  $i\in \ZZ$.
\end{dfn}

\begin{exa}
  \label{exa:dold}
  Every acyclic semi-flat complex is pure-acyclic, see
  \thmcite[7.3]{LWCHHl15}. On the other hand, the $\ZZ/4\ZZ$-complex
  known as the \emph{Dold complex},
  $$\cdots \lra \ZZ/4\ZZ \xra{2} \ZZ/4\ZZ \xra{2} \ZZ/4\ZZ \lra \cdots\:,$$
  is an acyclic complex of flat modules which is not pure
  acyclic. Indeed, the cycle submodules $2\ZZ/4\ZZ$ are torsion and
  hence not flat, cf.~\pgref{pure-m}.
\end{exa}

Recall that an $R$-module $P$ is \emph{pure-projective} if the
sequence
\begin{equation*}
  0 \lra \Hom{P}{L} \lra \Hom{P}{M} \lra \Hom{P}{N} \lra 0
\end{equation*}
is exact for every pure exact sequence $0 \to L \to M \to N \to 0$ of
$R$-modules.  Pure-injective modules are defined dually.

\begin{rmk}
  \label{rmk:emmanouil}
  An $R$-complex $M$ is by definition pure-acyclic if (and only if)
  $\Hom{A}{M}$ is acyclic for every finitely presented $R$-module
  $A$. Emmanouil \thmcite[3.6]{IEm16} shows that $M$ is pure-acyclic
  (if and) only if $\Hom{A}{M}$ is acyclic for every complex $A$ of
  pure-projective $R$-modules. Dually, $M$ is pure-acyclic if and only
  if $\Hom{M}{E}$ is acyclic for every pure-injective $R$-module $E$,
  equivalently, every complex $E$ of pure-injective $R$-modules. This
  was proved by {\v{S}}tov{\'{\i}}{\v{c}}ek \thmcite[5.4]{JSt}; see
  also Bazzoni, Cort\'es Izurdiaga, and Estrada \rmkcite[4.7]{BCE}.
\end{rmk}

From the proof of \corcite[2.6]{BCE} one can extract:

\begin{fct}
  \label{BG}
  Every pure-acyclic complex of pure-projective modules is
  contractible.
\end{fct}

For complexes of projective modules this follows from an earlier
result of Benson and Goodearl \thmcite[2.5]{DJBKRG00}; see
\prpcite[7.6]{LWCHHl15}. From the proof of \corcite[4.5]{BCE} one can
extract the dual result:

\begin{fct}
  \label{BCE}
  Every pure-acyclic complex of pure-injective modules is
  contractible.
\end{fct}

To close the section we apply these two facts to characterize von
Neumann regular rings in pure homological terms.

\begin{thm}
  \label{thm:vnr}
  The following conditions are equivalent.
  \begin{eqc}
  \item $R$ is von Neumann regular.
  \item Every acyclic $R$-complex is pure-acyclic.
  \item Every $R$-complex is semi-flat.
  \item Every complex of pure-projective $R$-modules is
    semi-projective.
  \item Every complex of pure-injective $R$-modules is semi-injective.
  \end{eqc}
\end{thm}

\begin{prf*}
  The (bi)implications $(i) \Leftrightarrow (ii)$ and
  $(iii) \Rightarrow (i)$ are clear from \rmkref{vnr}.

  $(i) \Rightarrow (iii)$: Let $M$ be an $R$-complex and $A$ be an
  acyclic $\Rop$-complex. As every $R$-module is flat, it follows from
  \corref{ptwize-C} that $\tp{A}{M}$ is acyclic, whence $M$ is
  semi-flat.
  
  $(ii) \Rightarrow (v)$: Let $J$ be a complex of pure-injective
  $R$-modules. By assumption every short exact sequence of $R$-modules
  is pure, so $J$ is a complex of injective modules. Let
  $\mapdef[\qra]{\iota}{J}{I}$ be a semi-injective resolution. The
  complex $\Cone{\iota}$ is an acyclic, hence pure-acyclic, complex of
  injective $R$-modules, so by \pgref{BCE} it is contractible. For an
  acyclic $R$-complex $A$, there is a triangle
  \begin{equation*}
    \Hom{A}{J} \lra \Hom{A}{I} \lra \Hom{A}{\Cone{\iota}} \lra
  \end{equation*}
  in the derived category of $\kk$-complexes. The middle complex is
  acyclic as $I$ is semi-injective, and the right-hand complexes is
  even contractible; it follows that $\Hom{A}{J}$ is acyclic, whence
  $J$ is semi-injective.

  $(v) \Rightarrow (i)$: A ring is von Neumann regular if and only if
  the opposite ring is so; it is thus sufficient to show that every
  $\Rop$-module is flat. The character dual $\Hom[\ZZ]{M}{\QZ}$ of an
  $\Rop$-module $M$ is a pure-injective $R$-module. It follows from
  the assumption that $\Hom[\ZZ]{M}{\QZ}$ is injective, whence $M$ is
  flat.

  $(iii) \Rightarrow (iv)$: Let $P$ be a complex of pure-projective
  $R$-modules. By assumption each module $P_i$ is flat and hence a
  pure quotient of a free $R$-module $L_i$, cf.~\pgref{pure-m}. By
  pure-projectivity of $P_i$, the homomorphism
  $\Hom{P_i}{L_i} \to \Hom{P_i}{P_i}$ is surjective, whence $P_i$ is a
  summand of $L_i$. Thus, $P$ is a complex of projective $R$-modules
  and semi-flat. Now proceed in parallel with the proof of
  $(ii) \Rightarrow (v)$ above, but invoke \pgref{BG} instead of
  \pgref{BCE}. Alternately see \thmcite[7.8]{LWCHHl15}.
  
  $(iv) \Rightarrow (ii)$: A finitely presented module is
  pure-projective, so by assumption finitely presented $R$-modules are
  projective. Now it follows from \rmkref{emmanouil} that every
  acyclic $R$-complex is pure-acyclic.
\end{prf*}

\section{Pure quasi-isomorphisms}
\label{sec:pqi}

\noindent
We continue with a discussion of morphisms with pure-acyclic mapping
cones.  Recall that a morphism $\gra$ of complexes is a homotopy
equivalence if and only if its mapping cone, $\Cone{\gra}$, is
contractible, while $\gra$ is a quasi-isomorphism if and only if
$\Cone{\gra}$ is acyclic. Pure quasi-isomorphisms are an intermediate
type of morphisms.

\begin{dfn}
  A morphism of $R$-complexes is called a \emph{pure
    quasi-isomorphism} if its mapping cone is a pure-acyclic complex.
\end{dfn}

\begin{rmk}
  \label{rmk:pqis}
  A morphism $\gra$ of $R$-complexes is a pure quasi-isomorphism if
  and only if $\Hom{A}{\gra}$ is a quasi-isomorphism for every
  finitely presented $R$-module $A$, equivalently for every complex
  $A$ of pure-projective $R$-modules. This follows in view of
  \rmkref{emmanouil} from the isomorphism
  $\Hom{A}{\Cone{\gra}} \is \Cone{\Hom{A}{\gra}}$.
\end{rmk}

\begin{exa}
  \label{exa:pqis}
  A contractible complex is pure-acyclic, so every homotopy
  equivalence is a pure quasi-isomorphism.  Further, every
  quasi-isomorphism of semi-flat complexes is a pure quasi-isomorphism
  by \corcite[7.4]{LWCHHl15} and \rmkref{pqis}.
\end{exa}

\begin{bfhpg}[Purity in the category of complexes]
  \label{pure-c}
  For $R$-complexes $L$ and $M$, the hom-set is denoted
  $\cathom{\cC(R)}{L}{M}$; it relates to the total Hom complex through
  the equality $\cathom{\cC(R)}{L}{M}=\Cy[0]{\Hom{L}{M}}$ of graded
  $\kk$-modules.

  An exact sequence of $R$-complexes $0 \to L \to M \to N \to 0$ is
  called \emph{pure} if the sequence of $\kk$-complexes
  \begin{equation*}
    0 \lra
    \cathom{\cC(R)}{A}{L} \lra \cathom{\cC(R)}{A}{M} \lra
    \cathom{\cC(R)}{A}{N} \lra 0
  \end{equation*}
  is exact for every bounded complex $A$ of finitely presented
  $R$-modules. It follows from \thmcite[5.1.3]{JGR99} and
  \thmcite[4.5]{LWCHHl15} that every pure exact sequence of
  $R$-complexes is degreewise pure.

  A subcomplex $L\subseteq M$, and a quotient complex $M/L$, are
  called pure if the canonical exact sequence
  $0 \to L \to M \to M/L \to 0$ is pure.
\end{bfhpg}

The next results help to recognize pure quasi-isomorphisms (homotopy
equivalences) and pure-acyclic pure subcomplexes (contractible split
subcomplexes).  In the literature, contractible complexes are at times
called split; we emphasize that by a (degreewise) split subcomplex we
mean a (degreewise) direct summand.

\begin{prp}
  \label{prp:ses_pa_vs_pqi}
  Let $0 \lra L\xra{\gra} M\xra{\grb} N \lra 0$ be a degreewise pure
  exact sequence in $\cC(R)$.  The following assertions hold.
  \begin{prt}
  \item The complex $L$ is pure-acyclic if and only if $\beta$ is a
    pure quasi-isomorphism.
  \item The complex $N$ is pure-acyclic if and only if $\alpha$ is a
    pure quasi-isomorphism.
  \end{prt}
  Moreover, if $L$ or $N$ is pure-acyclic, then the sequence is pure
  in $\cC(R)$.

  In particular, a pure-acyclic subcomplex is a pure subcomplex if and
  only if it is a degreewise pure subcomplex.
\end{prp}

\begin{prf*}
  For every finitely presented $R$-module $A$, the sequence
  \begin{equation*}
    0\lra \Hom{A}{L} \xrightarrow{\Hom{A}{\gra}} 
    \Hom{A}{M} \xrightarrow{\Hom{A}{\grb}} \Hom{A}{N}\lra 0
  \end{equation*}
  is exact, as the given exact sequence is degreewise pure. The
  complex $\Hom{A}{L}$ is acyclic if and only if $\Hom{A}{\grb}$ is a
  quasi-isomorphism, and $\Hom{A}{N}$ is acyclic if and only if
  $\Hom{A}{\gra}$ is a quasi-isomorphism. Now (a) and (b) follow from
  \rmkref{pqis}.
  
  Finally, given a bounded complex $A$ of finitely presented
  $R$-modules, we must verify exactness of the sequence
  \begin{equation*}
    0\lra \cathom{\cC(R)}{A}{L}\lra \cathom{\cC(R)}{A}{M}\lra 
    \cathom{\cC(R)}{A}{N}\lra 0\:.
  \end{equation*}
  As $\cathom{\cC(R)}{A}{-}$ is left exact and one has
  $\cathom{\cC(R)}{A}{-} = \Cy[0]{\Hom{A}{-}}$, this amounts to
  showing that the map
  \begin{equation*}
    \Cy[0]{\Hom{A}{M}}\xrightarrow{\Cy[0]{\Hom{A}{\beta}}} \Cy[0]{\Hom{A}{N}}
  \end{equation*}
  is surjective.  First notice that the morphism $\Hom{A}{\beta}$ is
  surjective as the given exact sequence is degreewise pure. If $L$ is
  pure-acyclic, then it follows from (a) and \rmkref{pqis} that
  $\Hom{A}{\beta}$ is a quasi-isomorphism and, therefore, surjective
  on cycles. If $N$ is pure-acyclic, then $\Hom{A}{N}$ is acyclic by
  \rmkref{emmanouil}. As every surjective chain map is surjective on
  boundaries, it follows that $\Hom{A}{\grb}$ is surjective on cycles.
\end{prf*}

\begin{prp}
  \label{prp:ses_c_vs_he}
  Let $0 \lra L\xra{\gra} M\xra{\grb} N \lra 0$ be a degreewise split
  exact sequence in $\cC(R)$. The following assertions hold.
  \begin{prt}
  \item The complex $L$ is contractible if and only if $\grb$ is a
    homotopy equivalence.
  \item The complex $N$ is contractible if and only if $\gra$ is a
    homotopy equivalence.
  \end{prt}
  Moreover, if $L$ or $N$ is contractible, then the sequence splits in
  $\cC(R)$.

  In particular, a contractible subcomplex is a split subcomplex if
  and only if it is a degreewise split subcomplex.
\end{prp}

\begin{prf*}
  Part (a) follows from \lemcite[1.6]{LLAAMr02}, and part (b) has an
  analogous proof.
\end{prf*}

In the sequel we use following two-of-three property of pure
quasi-isomorphisms.

\begin{lem}
  \label{lem:2of3pure}
  Let $\mapdef{\alpha}{L}{M}$ and $\mapdef{\beta}{M}{N}$ be morphisms
  of $R$-complexes. If any two of $\gra$, $\grb$, and $\grb\gra$ are
  pure quasi-isomorphisms, then so is the third.
\end{lem}

\begin{prf*}
  For every finitely presented $R$-module $A$ one has
  \begin{equation*}
    \H{\Hom{A}{\grb}}\H{\Hom{A}{\gra}} \deq \H{\Hom{A}{\grb\gra}}\:,
  \end{equation*}
  which shows that if any two of the morphisms $\H{\Hom{A}{\gra}}$,
  $\H{\Hom{A}{\grb}}$, and $\H{\Hom{A}{\grb\gra}}$ are isomorphisms,
  then so is the third. Now the statement follows from \rmkref{pqis}.
\end{prf*}

It is an elementary observation that an acyclic semi-projective
(-injective) complex is contractible and, hence, a quasi-isomorphism
of semi-projective (-injective) complexes is a homotopy
equivalence. In the same vein one has the following immediate
consequences of \pgref{BG} and \pgref{BCE}.

\begin{cor}
  \label{cor:pqi_of_proj}
  A pure quasi-isomorphism of complexes of pure-projective modules is
  a homotopy equivalence.\qed
\end{cor}

\begin{cor}
  \label{cor:pqi_of_inj}
  A pure quasi-isomorphism of complexes of pure-injective modules is a
  homotopy equivalence.\qed
\end{cor}

\section{Flavors of minimality}
\label{sec:flavors_of_min}

\noindent
We introduce the notion of \emph{pure-minimality} and explore how it
compares to notions found in the literature.  This section paves the
way for our main results in \secref{pm_replacements}.

We start by recalling that an $R$-complex $M$ is {\em minimal} if
every homotopy equivalence $M\to M$ is an isomorphism or,
equivalently, every morphism $M\to M$ that is homotopic to the
identity $1^M$ is an isomorphism.  In our context, this definition is
best known from \cite{LLAAMr02}. It is also an instance of Roig's
\cite{ARg93} notion of $S$-left or $S$-right minimality: the one where
$S$ in \defcite[1.1]{ARg93} is the class of homotopy equivalences.

Every complex has a minimal semi-injective resolution; see \cite{dga}
or \prpcite[B.2]{HKr05}. Minimal semi-projective resolutions are more
tricky: If $R$ is left perfect---such that flat $R$-modules are
projective---then every $R$-complex has a minimal semi-projective
resolution. If $R$ is semi-perfect---such that finitely generated flat
$R$-modules are projective---then every $R$-complex $M$ with $\H{M}$
degreewise finitely generated and $\H[i]{M} =0$ for $i \ll 0$ has a
minimal semi-projective resolution, see \cite{dga}. For the case of
resolutions of modules over a perfect ring, one can refer to Eilenberg
\cite{SEl56}.

\begin{exa}
  \label{exa:ZQ}
  The $\ZZ$-complex $F = 0 \lra \ZZ \xra{\iota} \QQ \lra 0$, where
  $\iota$ is the natural embedding, is minimal. Indeed, the morphisms
  $\pm1^F$ are the only homotopy equivalences $F \to F$, because there
  no nonzero homomorphisms $\QQ \to \ZZ$.
\end{exa}

The following weaker condition already detects minimality of complexes
of injective or projective modules; see
\prpref{implications}\footnote{\,The gist of this result is that
  minimality and split-minimality are equivalent notions for complexes
  of injective modules and complexes of projective modules. For the
  former this is true as stated, and for the latter it is true over
  rings where minimal complexes of projective modules are known to
  exist, in the strong sense that every complex of projectives
  decomposes as a direct sum of a minimal complex and a contractible
  one. The proof references the unpublished \cite{dga}, and the result
  serves to encapsulate this reference in the sense that one can
  henceforth focus on split-minimality of complexes of projective
  modules.}  below.

\begin{dfn}
  \label{dfn:split_min}
  An $R$-complex $M$ is called {\em split-minimal} if the zero complex
  is the only contractible split subcomplex of $M$.
\end{dfn}

\begin{exa}
  \label{exa:F}
  Let $p \ge 5$ be a prime and denote by $\ZZ_{p}$ the integers
  localized at the powers of $p$. The $\ZZ$-complex
  $\smash{F = 0 \lra \ZZ_{p} \xra{2} \ZZ_{p} \lra 0}$ is split-minimal
  but not minimal.  Indeed, every $\ZZ$-submodule of $\ZZ_{p}$ is
  cyclic, so any two non-zero submodules have a nonzero
  intersection. Hence $\ZZ_p$ has no non-trivial direct summands, and
  the zero complex is the only acyclic split subcomplex of $F$.
  However, the morphism $3^F$ is homotopic to $1^F$ (the homotopy is
  given by the identity on $\ZZ_p$) but not an isomorphism.
\end{exa}

To frame the next result we point out that the complex $F$ in
\exaref{F} is a complex of flat modules.

\begin{prp}
  \label{prp:implications}
  Let $M$ be an $R$-complex. If $M$ is minimal, then it is
  split-minimal; the converse holds under each of the following
  conditions:
  \begin{prt}
  \item $M$ is a complex injective $R$-modules.
  \item $R$ is left perfect; $M$ is a complex of projective
    $R$-modules.
  \item $R$ is semi-perfect; $M$ is a complex of finitely generated
    projective $R$-modules.
  \end{prt}
\end{prp}

\begin{prf*}  
  A minimal $R$-complex has by \prpcite[1.7]{LLAAMr02} no nonzero
  contractible split subcomplexes, hence it is split-minimal.
  
  (a): Let $M$ be a split-minimal complex of injective $R$-modules; it
  has by \cite{dga} or \prpcite[B.2]{HKr05} a decomposition
  $M= M'\oplus M''$, where $M'$ is minimal and $M''$ is
  contractible. Thus $M''$ is a contractible split subcomplex of $M$,
  and so $M''=0$.
  
  (b) \& (c): Suppose that $R$ is left perfect and $M$ is a
  split-minimal complex of projective $R$-modules, or that $R$ is
  semi-perfect and $M$ is a split-minimal complex of finitely
  generated projective $R$-modules. In either case, $M$ has a
  decomposition $M= M'\oplus M''$, where $M'$ is minimal and $M''$ is
  contractible; see \cite{dga}. As above, it follows that $M''=0$, and
  so $M$ is minimal.
\end{prf*}

In search of a useful notion of minimality for complexes of flat
modules, we turn to purity to formulate an analogue of
split-minimality. Given a semi-flat complex $F$ with flat dimension
$n$, the cokernels $\Co[i]{F}$ are flat for $i\ge n$, and it follows
that $P = \cdots \to F_{n+1}\to \Bo[n]{F}\to 0$ is a pure-acyclic
degreewise pure subcomplex of $F$.

At the least, a ``minimal'' semi-flat complex ought to vanish beyond
the flat dimension, but a minimal semi-flat complex need not meet that
requirement; see \exaref{minimal_not_puremin}. Inspired by this
example, we introduce a notion of minimality that forbids nonzero pure
acyclic pure (equivalently, degreewise pure) subcomplexes.

\begin{dfn}
  \label{dfn:M}
  An $R$-complex $M$ is called {\em pure-minimal} if the zero complex
  is the only pure-acyclic pure subcomplex of $M$.
\end{dfn}

\begin{rmk}
  \label{rmk:pm_wm}
  Every pure-minimal complex is split-minimal. Indeed every
  contractible complex is pure-acyclic and every split subcomplex is a
  pure subcomplex.
\end{rmk}

\begin{exa}
  \label{exa:ZQ'}
  The minimal $\ZZ$-complex $F$ from \exaref{ZQ} is pure-minimal. To
  see this, recall that a $\ZZ$-module is flat if and only if it is
  torsion-free; this means that $0$ and $\ZZ$ are the only pure
  submodules of $\ZZ$; cf.~\pgref{pure-m}. The only candidate for a
  nonzero pure-acyclic pure subcomplex of $F$ is, therefore,
  $0 \lra \ZZ \xra{=} \ZZ \lra 0$, but $\ZZ$ is not a pure submodule
  of $\QQ$, as $\QQ/\ZZ$ is torsion.
\end{exa}

\begin{exa}
  \label{exa:F'}
  The $\ZZ$-complex $F$ in \exaref{F} is pure-minimal but not
  minimal. Indeed, as every $\ZZ$-submodule of $\ZZ_p$ is cyclic, it
  follows that every non-trivial quotient is torsion.  As in
  \exaref{ZQ'} it then follows that the flat $\ZZ$-module $\ZZ_p$ has
  no non-trivial pure submodules. Thus, as $F$ is not acyclic, the
  zero complex is the only acyclic pure subcomplex of $F$.
\end{exa}

\begin{exa}
  \label{exa:minimal_not_puremin}
  Let $p$ be a prime and $\Zphat$ denote the $p$-adic completion of
  $\ZZ_{(p)}$.  The $\ZZ_{(p)}$-complex
  $F = 0\to \ZZ_{(p)} \to \Zphat \to 0$ is minimal but not
  pure-minimal.  To see this recall, for example from
  \lemcite[3.3]{AEL-02}, that one has
  $\Hom[\ZZ_{(p)}]{\Zphat}{\ZZ_{(p)}}=0$. It follows that every
  homotopy equivalence $F\to F$ is an isomorphism.  Thus, $F$ and
  $\Zphat/\ZZ_{(p)}$ are both minimal semi-flat resolutions of
  $\Zphat/\ZZ_{(p)}$; see \thmcite[(4.74)]{Lam2}. However, $F$ is not
  pure-minimal as the subcomplex
  $\smash{0\to \ZZ_{(p)} \xrightarrow{=} \ZZ_{(p)}\to 0}$ is
  degreewise pure and hence pure per \prpref{ses_pa_vs_pqi}.
\end{exa}

\begin{rmk}
  If one takes the class $S$ in \defcite[1.1]{ARg93} to be that of
  pure quasi-isomorphisms, then an $S$-right minimal complex in the
  sense of Roig is pure-minimal. Indeed, let $M$ be an $R$-complex and
  $P$ a pure-acyclic pure subcomplex of $M$; by \prpref{ses_pa_vs_pqi}
  the canonical map $\mapdef{\pi}{M}{M/P}$ is a pure
  quasi-isomorphism. Thus, if $M$ is $S$-right minimal, then $\pi$ has
  a left inverse, which implies that $\pi$ is injective and hence
  $P=0$.  The converse fails: with $F$ as in \exaref{F}, the map $3^F$
  is a pure quasi-isomorphism, see \exaref{pqis}, with no left
  inverse.
\end{rmk}

The next two corollaries paraphrase parts of
\prpref[Propositions~]{ses_pa_vs_pqi} and \prpref[]{ses_c_vs_he}.

\begin{cor}
  \label{cor:M1}
  Let $M$ be an $R$-complex. The next conditions are equivalent.
  \begin{eqc}
  \item $M$ is pure-minimal.
  \item The zero complex is the only pure-acyclic degreewise pure
    subcomplex of $M$.
  \item In a degreewise pure exact sequence
    $\smash{0 \lra L \lra M \xra{\grb} N \lra 0}$ the morphism $\grb$
    is a pure quasi-isomorphism if and only if it is an
    isomorphism. \qed
  \end{eqc}
\end{cor}

\begin{cor}
  \label{cor:M4}
  Let $M$ be an $R$-complex. The next conditions are equivalent.
  \begin{eqc}
  \item $M$ is split-minimal.
  \item The zero complex is the only contractible degreewise split
    subcomplex of $M$.
  \item In a degreewise split exact sequence
    $\smash{0 \lra L \lra M \xra{\grb} N\lra 0}$ the morphism $\grb$
    is a homotopy equivalence if and only if it is an isomorphism.\qed
  \end{eqc}
\end{cor}




We now show that split-minimality and pure-minimality coincide in
standard settings while we already saw in \exaref{minimal_not_puremin}
that a (split-)minimal complex of flat modules need not be
pure-minimal.

\begin{thm}
  \label{thm:pureminiffweakmin}
  Let $M$ be an $R$-complex. Under each of the conditions \mbox{{\rm
      (a)--(d)}} below, every pure-acyclic pure subcomplex of $M$ is
  contractible and a split subcomplex; in particular, the complex $M$
  is split-minimal if and only if it is pure-minimal.
  \begin{prt}
  \item $R$ is left noetherian; $M$ is a complex of injective
    $R$-modules.
  \item $R$ is left perfect; $M$ is a complex of projective
    $R$-modules.
  \item $R$ is semi-perfect; $M$ is a complex of finitely generated
    projective $R$-modules.
  \item $R$ is left noetherian; $M$ is a complex of finitely generated
    projective~\mbox{$R$-modules}.
  \end{prt}
\end{thm}

\begin{prf*}
  Let $P$ be a pure-acyclic pure subcomplex of $M$ and consider the
  pure exact sequence
  \begin{equation*}
    \tag{$\ast$}
    0 \lra P \lra M \to M/P \lra 0\:.
  \end{equation*}
  Under any one of the assumptions (a)--(d) the sequence is degreewise
  split exact by \lemref{pmsm}. Applied to the exact sequences
  $0 \to \Cy[i]{P} \to P_i \to \Cy[i-1]{P} \to 0$ the same lemma shows
  that $P$ is contractible. Now it follows from \prpref{ses_c_vs_he}
  that the sequence $(\ast)$ splits.

  By \rmkref{pm_wm} every pure-minimal complex is split-minimal, and
  the argument above shows that the converse holds under each of the
  conditions (a)--(d).
\end{prf*}

\begin{cor}
  \label{cor:standardsituations}
  Let $R$ and $M$ be as in {\rm \thmref[]{pureminiffweakmin}(a), (b),}
  or {\rm (c)}. The $R$-complex $M$ is pure-minimal if and only if it
  is split-minimal if and only if it is minimal.
\end{cor}

\begin{prf*}
  Combine \prpref{implications} and \thmref{pureminiffweakmin}.
\end{prf*}

\begin{ipg}
  The diagram below summarizes the (non-)implications among the
  notions of minimality considered in this section. We stress that
  while the three notions agree under each of the assumptions (a),
  (b), or (c) in \thmref{pureminiffweakmin}, the examples that lie
  behind the non-implications deal with semi-flat complexes over PIDs.
  
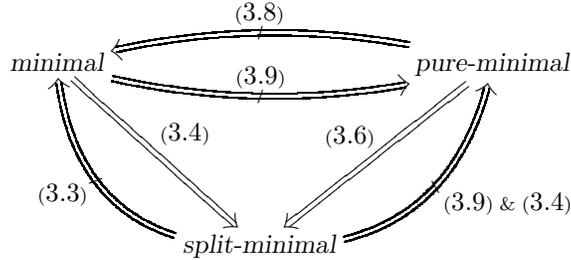
\begin{equation*}
    \xymatrix{
      \textsl{minimal} 
      \ar@/_1pc/@{=>}[rr]|-{\SelectTips{cm}{}\object@{/}}|{}^-{\pgref{exa:minimal_not_puremin}} 
      \ar@{=>}[ddr]^-{\pgref{prp:implications}} && 
      \textsl{pure-minimal}
      \ar@{=>}[ddl]_-{\pgref{rmk:pm_wm}} 
      \ar@/_1pc/@{=>}[ll]|-{\SelectTips{cm}{}\object@{/}}|{}_-{\pgref{exa:F'}}\\ &&\\
      &\textsl{split-minimal} 
      \ar@/^2pc/@{=>}[uul]|-{\SelectTips{cm}{}\object@{/}}|{}^-{\pgref{exa:F}} 
      \ar@/_2pc/@{=>}[uur]|-{\SelectTips{cm}{}
        \object@{/}}|{}_-{\pgref{exa:minimal_not_puremin} \;\&\; \pgref{prp:implications}}
    }
  \end{equation*}
\end{ipg}

\section{Minimality of acyclic complexes}\label{sec:min_ac}

\noindent
The zero complex is minimal as can be; this short section complements
the preceding one by spelling out what the flavors of minimality mean
for acyclic complexes.

We start by noticing that the Dold complex from \exaref{dold} is an
acyclic complex of projective and injective modules which is both
minimal and pure-minimal, as $\ZZ/4\ZZ$ has no non-trivial pure
submodule. Thus, nonzero minimal and pure-minimal acyclic complexes
exist over quasi-Frobenius rings.

\begin{prp}
  \label{prp:asm}
  The zero complex is the only contractible split-minimal complex. In
  particular, the zero complex is the only
  \begin{itemize}
  \item acyclic split-minimal semi-injective complex,
  \item acyclic split-minimal semi-projective complex,
  \item pure-acyclic split-minimal complex of pure-injective modules,
  \item pure-acyclic split-minimal complex of pure-projective modules.
  \end{itemize}
\end{prp}

\begin{prf*}
  The first assertion is immediate from \dfnref{split_min} and the
  remaining follow from \cite[1.3.P and 1.3.I]{LLAHBF91}, \pgref{BG},
  and \pgref{BCE}.
\end{prf*}

\begin{exa}
  Assume that $R$ is semi-simple; that is, every acyclic $R$-complex
  is contractible.  It follows from \prpref{asm} that the zero complex
  is the only acyclic split-minimal $R$-complex. Even more, an
  $R$-complex $M$ is split-minimal if and only if the zero complex is
  the only acyclic subcomplex of $M$. Indeed, the ``if'' is trivial
  and the ``only if'' follows from \prpref{ses_c_vs_he}.
\end{exa}

\begin{prp}
  \label{prp:apm}
  The zero complex is the only pure-acyclic pure-minimal complex; in
  particular it is the only acyclic pure-minimal semi-flat complex.
\end{prp}

\begin{prf*}
  The first assertion is immediate from \dfnref{M} and the second
  follows from \thmcite[7.3]{LWCHHl15}.
\end{prf*}

\begin{exa}
  \label{exa:vnr}
  Assume that $R$ is von Neumann regular; that is, every acyclic
  $R$-complex is pure-acyclic; see \thmref{vnr}. It follows from
  \prpref{apm} that the zero complex is the only acyclic pure-minimal
  $R$-complex. (In fact, this property characterizes von Neumann
  regular rings; see \corref{vnr}.) Even more, an $R$-complex $M$ is
  pure-minimal if and only if the zero complex is the only acyclic
  subcomplex of $M$. Here ``only if'' follows from
  \prpref{ses_pa_vs_pqi} and ``if'' is clear.
\end{exa}

For work in the derived category of chain complexes---computation of
derived functors for example---the emphasis is on distinguished
complexes of injective (projective or flat) modules, namely the
semi-injective (-projective or -flat) complexes.  The next result
shows that the noetherian hypothesis in \thmref{pureminiffweakmin}(a),
so to speak, does not impact work in the derived category.

\begin{prp}
  \label{prp:puremin-semiinj}
  Let $I$ be a semi-injective $R$-complex. The following conditions
  are equivalent:
  \begin{eqc}
  \item $I$ is minimal.
  \item $I$ is split-minimal.
  \item $I$ is pure-minimal.
  \item The zero complex is the only acyclic subcomplex of $I$.
  \end{eqc}
\end{prp}

\begin{prf*}
  In view of \prpref{implications}(a) and \rmkref{pm_wm} it suffices
  to show that $(i)$ implies $(iv)$. Assume that $I$ is minimal and
  let $A$ be an acyclic subcomplex. As $I$ is semi-injective, the
  quasi-isomorphism $\mapdef{\pi}{I}{I/A}$ has a left inverse up to
  homotopy. That is, there is a morphism $\mapdef{\gamma}{I/A}{I}$
  such that $\gamma\pi$ is homotopic to $1^I$. As $I$ is minimal, it
  follows that $\gamma\pi$ is an isomorphism and, therefore, $A=0$.
\end{prf*}

Pure-minimal semi-flat complexes have a similar
characterization.\footnote{\,When they exist, minimal semi-projective
  complexes are characterized by having only the trivial acyclic
  quotient complex.}

\begin{prp}
  \label{prp:puremin-semiflat}
  Let $F$ be a semi-flat $R$-complex. The following conditions are
  equivalent:
  \begin{eqc}
  \item $F$ is pure-minimal.
  \item The zero complex is the only acyclic pure subcomplex of $F$.
  \end{eqc}
\end{prp}

\begin{prf*}
  It is clear that $(ii)$ implies $(i)$. For the converse, let $P$ be
  an acyclic pure subcomplex of $F$. It follows from
  \prpcite[6.2]{LWCHHl15} that $P$ is semi-flat, hence $P$ is pure
  acyclic by \thmcite[7.3]{LWCHHl15} and, therefore, $P=0$.
\end{prf*}

\section{Pure-minimal replacements}
\label{sec:pm_replacements}

\noindent
Every chain complex $M$ has a minimal, hence pure-minimal,
semi-injective resolution $M \qra I$; see \prpref{puremin-semiinj}. In
particular, $M$ and $I$ are isomorphic in the derived category. We
proceed to show that every complex is isomorphic, in the derived
category, to a pure-minimal semi-flat complex.

The gist of the next theorem, which is our central construction, is
that every complex has a pure-acyclic pure subcomplex, such that the
associated quotient is pure-minimal.

\begin{thm}
  \label{thm:sub_quotient_complexes}
  Let $M$ be an $R$-complex. There is a pure exact sequence in
  $\cC(R)$
  \begin{equation*}
    0\lra P \lra M \lra M/P \lra 0
  \end{equation*}
  with $P$ pure-acyclic and $M/P$ pure-minimal. Consequently, the map
  $M\lra M/P$ is a pure quasi-isomorphism.
\end{thm}

\begin{prf*}
  Consider the set of all pure-acyclic pure subcomplexes of $M$,
  ordered by containment. Let $\Lambda$ be a chain in this set, and
  let $U$ be the union $\colim_{A\in \Lambda} A$.  It is standard that
  $U$ is a subcomplex of $M$ and we proceed to show that it is a pure
  subcomplex and pure-acyclic.

  First we verify that $U$ is a pure subcomplex. Let $F$ be a bounded
  complex of finitely presented $R$-modules. For every $A\in \Lambda$
  there is an exact sequence
  \begin{equation*}
    0\lra \cathom{\cC(R)}{F}{A}\lra \cathom{\cC(R)}{F}{M}\lra 
    \cathom{\cC(R)}{F}{M/A}\lra 0\:.
  \end{equation*}
  Recall, e.g.\ from \thmcite[4.5]{LWCHHl15}, that as $\Lambda$ is
  filtered there is a natural isomorphism,
  \begin{equation*}
    \colim_{A\in \Lambda}\cathom{\cC(R)}{F}{-}\ira
    \cathom{\cC(R)}{F}{\colim_{A\in \Lambda}(-)}\:.
  \end{equation*}
  Using this, along with the fact that filtered colimits are exact, we
  obtain the following commutative diagram with exact rows:
  \begin{equation*}
    \xymatrix@C=.825pc{
      0 \ar[r]
      & \smash{\displaystyle\colim_{A\in \Lambda}} \cathom{\cC(R)}{F}{A} \ar[r] \ar[d]^-\is
      & \smash{\displaystyle\colim_{A\in \Lambda}} \cathom{\cC(R)}{F}{M} \ar[r] \ar[d]^-\is
      & \smash{\displaystyle\colim_{A\in \Lambda}} \cathom{\cC(R)}{F}{M/A} \ar[r] \ar[d]^-\is
      & 0\\
      0 \ar[r]
      & \cathom{\cC(R)}{F}{\displaystyle\colim_{A\in \Lambda}A} \ar[r]
      & \cathom{\cC(R)}{F}{\displaystyle\colim_{A\in \Lambda}M} \ar[r]
      & \cathom{\cC(R)}{F}{\displaystyle\colim_{A\in \Lambda}M/A}
      & 
      \\
    }
  \end{equation*}
  It follows from a simple diagram chase that
  \begin{equation*}
    \cathom{\cC(R)}{F}{\displaystyle\colim_{A\in \Lambda}M} \lra 
    \cathom{\cC(R)}{F}{\displaystyle\colim_{A\in \Lambda}M/A}\lra 0
  \end{equation*}
  is exact. In view of the canonical isomorphisms
  \begin{equation*}
    M/\colim_{A\in \Lambda}A \dis 
    \colim_{A\in \Lambda}M/\colim_{A\in \Lambda}A
    \dis \colim_{A\in \Lambda}M/A 
  \end{equation*}
  it thus follows that $U=\colim_{A\in \Lambda}A$ is a pure subcomplex
  of $M$.
  
  Next we argue that $U$ is pure-acyclic. Let $F$ be a finitely
  presented $R$-module; we have to show that $\Hom{F}{U}$ is
  acyclic. Since the functor $\Hom{F}{-}$ preserves filtered colimits,
  one has
  \begin{equation*}
    \Hom{F}{U} \deq \Hom{F}{\colim_{A\in \Lambda}A} \dis 
    \colim_{A\in \Lambda}\Hom{F}{A}\:.
  \end{equation*}
  As each complex $A$ is pure-acyclic, the complexes $\Hom{F}{A}$ are
  acyclic. Finally, $\colim_{A\in \Lambda}(-)$ preserves acyclicity
  and so $\Hom{F}{U}$ is acyclic.
  
  By Zorn's lemma, there exists a maximal pure-acyclic pure subcomplex
  $P$ of~$M$. To show that $M/P$ is pure-minimal, let
  $P'/P \subseteq M/P$ be a pure-acyclic pure subcomplex. Consider the
  commutative diagram with exact rows and columns:
  \begin{equation*}
    \xymatrix@=1.5pc{
      & 0 \ar[d] & 0 \ar[d]\\
      & P \ar[d] \ar@{=}[r]& P \ar[d] \\
      0 \ar[r] & P' \ar[r]\ar[d] & M \ar[r]\ar[d] & M/P' \ar@{=}[d]\ar[r] & 0\\
      0 \ar[r] & P'/P \ar[r] \ar[d] & M/P \ar[r]\ar[d] & M/P' \ar[r] & 0\\
      & 0 & 0
    }
  \end{equation*}
  The bottom row and the middle column are pure sequences in $\cC(R)$
  by the assumptions and by what we have shown above. Hence, for every
  bounded complex $F$ of finitely presented $R$-modules, application
  of $\cathom{\cC(R)}{F}{-}$ yields another commutative diagram with
  exact rows and columns
  \begin{equation*}
    \xymatrix@=1.5pc{
      & 0 \ar[d] & 0 \ar[d]\\
      & \cathom{\cC(R)}{F}{P} \ar[d] \ar@{=}[r]
      & \cathom{\cC(R)}{F}{P} \ar[d] \\
      0 \ar[r] & \cathom{\cC(R)}{F}{P'} \ar[r]\ar[d] 
      & \cathom{\cC(R)}{F}{M} \ar[r]\ar[d] 
      & \cathom{\cC(R)}{F}{M/P'} \ar@{=}[d] \\
      0 \ar[r] & \cathom{\cC(R)}{F}{P'/P} \ar[r] 
      & \cathom{\cC(R)}{F}{M/P} \ar[r]\ar[d] 
      & \cathom{\cC(R)}{F}{M/P'} \ar[r] & 0\\
      & & 0
    }
  \end{equation*}
  A diagram chase shows that the morphism
  $\cathom{\cC(R)}{F}{M} \lra \cathom{\cC(R)}{F}{M/P'}$ is surjective,
  whence $P'$ is a pure subcomplex of $M$. To see that $P'$ is
  pure-acyclic it is now by \prpref{ses_pa_vs_pqi} sufficient to show
  that the canonical map $M \to M/P'$ is a pure
  quasi-isomorphism. This map is the composite of canonical maps
  \begin{equation*}
    M \lra M/P \lra M/P'\:,
  \end{equation*}
  both of which are pure quasi-isomorphisms, again by
  \prpref{ses_pa_vs_pqi}. Now it follows from \lemref{2of3pure} that
  $M \to M/P'$ is a pure quasi-isomorphism. As $P$ is a maximal pure
  acyclic pure subcomplex of $M$, one gets $P'/P=0$, and it follows
  that $M/P$ is pure-minimal.
  
  An application of \prpref{ses_pa_vs_pqi} now shows that the map
  $M\to M/P$ is a pure quasi-isomorphism.
\end{prf*}

\begin{cor}
  \label{cor:vnr}
  The following conditions are equivalent.
  \begin{eqc}
  \item $R$ is von Neumann regular.
  \item The zero complex is the only acyclic pure-minimal $R$-complex.
  \item An $R$-complex $M$ is pure-minimal if and only if the zero
    complex is the only acyclic subcomplex of $M$.
  \end{eqc}
\end{cor}

\begin{prf*}
  It was noted in \exaref{vnr} that $(iii)$ follows from $(i)$ by way
  of \prpref{ses_pa_vs_pqi}, and $(iii)$ clearly implies $(ii)$. To
  see that $(ii)$ implies $(i)$, let $M$ be an acyclic $R$-complex. By
  \thmref{sub_quotient_complexes} there is a pure-acyclic subcomplex
  $P$ of $M$ such that the quotient $M/P$ is pure-minimal and
  acyclic. By assumption, $M/P=0$ so $M=P$ is pure-acyclic. It now
  follows from \thmref{vnr} that $R$ is von Neumann regular.
\end{prf*}

In classical settings, such as in \thmref{pureminiffweakmin}(a,b,c), a
complex decomposes as a direct sum of a minimal complex and a
contractible one. These decompositions are recovered by
\thmref{sub_quotient_complexes}---together with
\corref{standardsituations}---which also yields a similar
decomposition for complexes of finitely generated projective modules
over noetherian rings.



\begin{cor}
  \label{cor:decomposition}
  Let $R$ and $M$ be as in {\rm \thmref[]{pureminiffweakmin}(a), (b),
    (c)}, or {\rm (d)}. The exact sequence $0\to P\to M\to M/P\to 0$
  from {\rm (\ref{thm:sub_quotient_complexes})} is split in $\cC(R)$
  and yields a decomposition $M\is P\oplus (M/P)$, where $P$ is
  contractible and $M/P$ is pure-minimal.
\end{cor}

\begin{prf*}
  Immediate from \thmref{pureminiffweakmin}.
\end{prf*}

\begin{thm}
  \label{thm:pm_pd}
  Let $R$ be left noetherian and $M$ be an $R$-complex with $\H{M}$
  degreewise finitely generated and $\H[i]{M}=0$ for $i \ll 0$. There
  is a semi-projective resolution $L \qra M$ with $L$ pure-minimal and
  degreewise finitely generated.  Furthermore, for every such
  resolution $L \qra M$ one~has
  \begin{equation*}
    \pd{M} \deq \sup\setof{i}{L_i\not=0}\:.
  \end{equation*}
\end{thm}

\begin{prf*}
  Notice first that if $M$ is acyclic, then one can take $L=0$.
  Assume that $M$ is not acyclic, and let $L' \qra M$ be a
  semi-projective resolution with $L'$ degreewise finitely generated;
  see \rmkcite[1.7]{LLAHBF91}. By \corref{decomposition} the complex
  $L'$ has a pure-minimal summand $L$. The complex $L$ is
  semi-projective, degreewise finitely generated, and isomorphic to
  $M$ in the derived category. By \cite[1.4.P]{LLAHBF91} there is a
  quasi-isomorphism $L \to M$.
  
  Let $L \qra M$ be a semi-projective resolution with $L$ pure-minimal
  and degreewise finitely generated. By \thmcite[2.4.P]{LLAHBF91} one
  has $\pd{M}\le \sup\setof{i}{L_i\not=0}$, and equality holds
  trivially if $\pd{M}=\infty$. If $M$ has finite projective dimension
  $n$, then the complex
  $L_{\subseteq n}=0\to \Co[n]{L}\to L_{n-1}\to \cdots$ is by
  \thmcite[2.4.P]{LLAHBF91} a split subcomplex of $L$ with
  contractible complement. It follows that the complement is zero,
  whence one has $L= L_{\subseteq n}$ and
  $\pd{M}=n=\sup\setof{i}{L_i\not=0}$.
\end{prf*}

\begin{rmk}
  For any pure-minimal degreewise finitely generated semi-projective
  complex $P$, the proof of \thmref{pm_pd} yields
  $\pd{P}=\sup\setof{i}{P_i\not=0}$.
\end{rmk}

The construction in \thmref{sub_quotient_complexes} also applies to
yield a complex that detects flat dimension. We recall from
\exacite[2.9.F]{LLAHBF91} that a complex $M$ need not have a semi-flat
resolution that detects its flat dimension, hence we settle for a {\em
  semi-flat replacement} of $M$, i.e.\ a semi-flat complex isomorphic
to $M$ in the derived category.

\begin{thm}
  \label{thm:pm_fd}
  For every $R$-complex $M$ there exists a pure-minimal semi-flat
  $R$-complex $F$ isomorphic to $M$ in the derived
  category. Furthermore, for every such complex $F$ one~has
  \begin{equation*}
    \fd{M} \deq \sup\setof{i}{F_i\not=0}\:.
  \end{equation*}
\end{thm}

\begin{prf*}
  If $M$ is acyclic, then one can take $F=0$, so assume that $M$ is
  not acyclic.  Let $L \qra M$ be a semi-projective
  resolution. \thmref{sub_quotient_complexes} yields a pure-acyclic
  pure subcomplex $P$ of $L$ such that the quotient $F=L/P$ is
  pure-minimal. As $L$ is semi-flat, it follows from
  \prpcite[6.2]{LWCHHl15} that $F$ is semi-flat as well. There are now
  quasi-isomorphisms
  \begin{equation*}
    M \qla L \qra F\:,
  \end{equation*}
  so $M$ and $F$ are isomorphic in the derived category.
  
  Let $F$ be a pure-minimal semi-flat $R$-complex isomorphic to $M$ in
  the derived category. By \thmcite[2.4.F]{LLAHBF91} one has
  $\fd{M}\le \sup\setof{i}{F_i\not=0}$, and equality holds trivially
  if $\fd{M}=\infty$. If $M$ has finite flat dimension $n$, then
  $F_{\subseteq n}=0\to \Co[n]{F}\to F_{n-1}\to \cdots$ is a semi-flat
  $R$-complex isomorphic to $M$ in the derived category; see
  \thmcite[2.4.F]{LLAHBF91}. Set $K = \Ker(F\to F_{\subseteq n})$ and
  consider the exact sequence
  \begin{equation*}
    0\lra K\lra F\lra F_{\subseteq n}\lra 0\:.
  \end{equation*}
  It is degreewise pure as $F_{\subseteq n}$ is a complex of flat
  modules. The morphism $F\lra F_{\subseteq n}$ is a pure
  quasi-isomorphism, see \exaref{pqis}, so it follows from
  \prpref{ses_pa_vs_pqi} that $K$ is a pure-acyclic pure subcomplex of
  $F$.  Since $F$ is pure-minimal, this means $K=0$. Hence one has
  $F= F_{\subseteq n}$ and $\fd{M}=n=\sup\setof{i}{F_i\not=0}$.
\end{prf*}

Minimal semi-projective resolutions are only known to exist for all
$R$-complexes if $R$ is a left perfect ring. We close this section
with a characterization of such rings in terms of existence of
pure-minimal semi-projective resolutions.

\begin{thm}
  \label{thm:perfect_pm}
  The following conditions on $R$ are equivalent.
  \begin{eqc}
  \item $R$ is left perfect.
  \item Every semi-flat $R$-complex is semi-projective.
  \item Every $R$-complex has a pure-minimal semi-projective
    resolution.
  \end{eqc}
\end{thm}

\begin{prf*}
  Every flat module over a perfect ring is projective, and a semi-flat
  complex of projective modules is semi-projective; see
  \thmcite[7.8]{LWCHHl15}. Thus $(i)$ implies $(ii)$. By
  \thmref{pm_fd} every $R$-complex $M$ has a pure-minimal semi-flat
  replacement $F$. Assuming $(ii)$ the complex $F$ is semi-projective,
  and it follows from \cite[1.4.P]{LLAHBF91} that there is a
  quasi-isomorphism $F \to M$. Thus $(ii)$ implies $(iii)$. To finish
  the proof, let $F$ be a flat $R$-module with pure-minimal
  semi-projective resolution $\mapdef[\qra]{\pi}{P}{F}$. As
  $\mapdef{\H{\pi}}{\H{P}}{\H{F}=F}$ is an isomorphism, $\pi$ is
  surjective, and it follows from \pgref{pure-m} that $K = \Ker{\pi}$
  is a degreewise pure subcomplex of $P$. Since $P$ and $F$ are
  semi-flat complexes, $\pi$ is a pure quasi-isomorphism, see
  \exaref{pqis}. From \corref{M1} it now follows that $\pi$ is an
  isomorphism. Thus every flat $R$-module is projective, whence $R$ is
  left perfect.
\end{prf*}

\appendix
\section*{Appendix. Sufficient conditions for acyclicity}
\stepcounter{section}

\noindent We collect a few technical results that are useful for
establishing acyclicity of Hom and tensor product complexes. The results
complement and improve those in \seccite[2]{CFH-06}; the proofs extend
and dualize an argument by Emmanouil \cite{IEm16}.

\begin{prp}
  \label{prp:ptwize-C}
  Let $M$ and $N$ be $R$-complexes. The complex $\Hom{M}{N}$ is
  acyclic if the following conditions are satisfied.
  \begin{prt}
  \item $\Hom{M_i}{N}$ is acyclic for every $i\in \ZZ$, and
  \item $\Hom{\Co[i]{M}}{N}$ is acyclic for every $i \ll 0$.
  \end{prt}
\end{prp}

\begin{prf*}
  Emmanouil's argument for \lemcite[2.6]{IEm16} can be adapted to
  apply; see also the argument for the dual result \prpref{ptwize-Z}
  below.
\end{prf*}

\begin{cor}
  \label{cor:ptwize-C}
  Let $L$ be an $\Rop$-complex and $M$ be an $R$-complex. The complex
  $\tp{L}{M}$ is acyclic if the next conditions are~satisfied.
  \begin{prt}
  \item $\tp{L}{M_i}$ is acyclic for every $i\in \ZZ$, and
  \item $\tp{L}{\Co[i]{M}}$ is acyclic for every $i \ll 0$.
  \end{prt}
\end{cor}

\begin{prf*}
  Recall that the complex $\tp{L}{M}$ is acyclic if and only if the
  dual complex $\Hom[\ZZ]{\tp{L}{M}}{\QQ/\ZZ}$ is acyclic. The result
  now follows from \prpref{ptwize-C} by way of Hom--tensor adjunction.
\end{prf*}

\begin{prp}
  \label{prp:ptwize-Z}
  Let $M$ and $N$ be $R$-complexes. The complex $\Hom{M}{N}$ is
  acyclic if the following conditions are satisfied.
  \begin{prt}
  \item $\Hom{M}{N_i}$ is acyclic for every $i\in \ZZ$, and
  \item $\Hom{M}{\Cy[i]{N}}$ is acyclic for every $i \gg 0$.
  \end{prt}
\end{prp}

\begin{prf*}
  It is well-known, for example from \lemcite[2.5]{CFH-06}, that
  condition (a) implies that the complex $\Hom{M}{N_{< m}}$ is acyclic
  for every $m$. From an application of $\Hom{M}{-}$ to the degreewise
  split exact sequence $0 \to N_{< m}\to N \to N_{\ge m} \to 0$ it
  follows that it is sufficient to prove that $\Hom{M}{N_{\ge m}}$ is
  acyclic for some integer $m$. Thus, without loss of generality
  assume that $N_i=0$ holds for $i \ll 0$ and that
  $\Hom{M}{\Cy[i]{N}}$ is acyclic for every $i \in \ZZ$.

  A homomorphism $M\to N$ is a cycle in $\Hom{M}{N}$ if and only if it
  is a chain map and a boundary if and only if it is
  null-homotopic. Let $\mapdef{\f}{M}{N}$ be a chain map; after
  shifting and reindexing we may assume that $\f$ has degree zero. The
  goal is to construct a homotopy from $\f$ to $0$, i.e.\ a family of
  $R$-module homomorphisms $\mapdef{\s_i}{M_i}{N_{i+1}}$ with
  $\f_i = \dif[i+1]{N}\s_i + \s_{i-1}\dif[i]{M}$. Evidently $\s_i$ has
  to be zero for $i \ll 0$; this provides the basis for an induction
  argument. Fix $n$ and assume that the desired homomorphisms $\s_i$
  have been constructed for $i \le n-2$; assume further that a
  homomorphism $\mapdef{\tau_{n-1}}{M_{n-1}}{N_n}$ with
  $\f_{n-1} = \dif[n]{N}\tau_{n-1} + \s_{n-2}\dif[n-1]{M}$ has been
  constructed. The map $\tau_{n-1}$ may not have all the properties
  required of $\s_{n-1}$, but in the induction step it is modified to
  yield the desired $\s_{n-1}$. For $i \ll 0$ one takes
  $\tau_i=0$. The next diagram depicts the data from the induction
  hypothesis.
  \begin{equation*}
    \xymatrix{
      \cdots \ar[r] & M_{n+1} \ar[r] \ar[d]^-{\f_{n+1}} 
      & M_{n} \ar[r] \ar[d]^-{\f_{n}} 
      & M_{n-1} \ar[r] \ar[d]^-{\f_{n-1}} \ar[dl]^-{\tau_{n-1}} 
      & M_{n-2} \ar[r] \ar[d]^-{\f_{n-2}} \ar[dl]^-{\s_{n-2}}
      & M_{n-3} \ar[r] \ar[d]^-{\f_{n-3}} \ar[dl]^-{\s_{n-3}} & \cdots\\
      \cdots \ar[r] & N_{n+1} \ar[r] & N_{n} \ar[r] & N_{n-1} \ar[r] 
      & N_{n-2} \ar[r] & N_{n-3} \ar[r] & \cdots
    }
  \end{equation*}
  In the induction step we need to construct homomorphisms
  \begin{equation*}
    \tag{1}
    \dmapdef{\tau_{n}}{M_{n}}{N_{n+1}}   \qand \dmapdef{\upsilon_{n-1}}{M_{n-1}}{N_n}
  \end{equation*}
  such that $\tau_n$ and $\s_{n-1} = \tau_{n-1} + \upsilon_{n-1}$ satisfy
  \begin{equation*}
    \tag{2}
    \f_{n} \deq \dif[n+1]{N}\tau_{n} + \s_{n-1}\dif[n]{M} \qand     
    \f_{n-1} \deq \dif[n]{N}\s_{n-1} + \s_{n-2}\dif[n-1]{M}\:.
  \end{equation*}
  In the next computation, the first equality holds as $\f$ is a chain
  map, and the second follows from the assumption on $\tau_{n-1}$,
  \begin{equation*}
    \dif[n]{N}(\f_n - \tau_{n-1}\dif[n]{M}) 
    \deq (\f_{n-1} - \dif[n]{N}\tau_{n-1})\dif[n]{M} 
    \deq \s_{n-2}\dif[n-1]{M}\dif[n]{M} \deq 0\:.
  \end{equation*}
  This shows that $\f_n - \tau_{n-1}\dif[n]{M}$ corestricts to a
  homomorphism $M_n \to \Cy[n]{N}$.  It is elementary to verify that
  the diagram
  \begin{equation*}
    \xymatrix{
      \cdots \ar[r] 
      & M_{n+3} \ar[r] \ar[d] 
      & M_{n+2} \ar[r] \ar[d]^-{\f_{n+1}\dif[n+2]{M}} 
      & M_{n+1} \ar[r] \ar[d]^-{\f_{n+1}}
      & M_{n} \ar[r] \ar[d]^-{\f_{n} -\tau_{n-1}\dif[n]{M}}
      & M_{n-1} \ar[r] \ar[d] & \cdots\\
      & 0 \ar[r] & \Cy[n+1]{N} \ \ar@{>->}[r] & N_{n+1} \ar[r] 
      & \Cy[n]{N} \ar[r] & 0
    }
  \end{equation*}
  is commutative; that is, the vertical maps form a chain map
  $\mapdef{\f'}{M}{N'}$, where $N'$ denotes the bottom row in the diagram. By the assumptions and \lemcite[2.5]{CFH-06}
  the complex $\Hom{M}{N'}$ is acyclic, so $\f'$ is null-homotopic. In
  particular, there exist homomorphisms $\tau_n$ and $\upsilon_{n-1}$
  as in $(1)$ with
  \begin{equation*}
    \f_{n} -\tau_{n-1}\dif[n]{M} \deq
    \dif[n+1]{N}\tau_n + \upsilon_{n-1}\dif[n]{M}\:.
  \end{equation*}
  It is now straightforward to verify that the identities in $(2)$
  hold.
\end{prf*}

\begin{cor}
  \label{cor:ptwize-Z}
  Let $L$ be an $\Rop$-complex and $M$ be an $R$-complex. The complex
  $\tp{L}{M}$ is acyclic if the next conditions are~satisfied.
  \begin{prt}
  \item $\tp{L_i}{M}$ is acyclic for every $i\in \ZZ$, and
  \item $\tp{\Bo[i]{L}}{M}$ is acyclic for every $i \ll 0$.
  \end{prt}
\end{cor}

\begin{prf*}
  As in the proof of \corref{ptwize-C} it suffices to show that the
  complex
  \begin{equation*}
    \Hom[\ZZ]{\tp{L}{M}}{\QQ/\ZZ} \dis \Hom{M}{\Hom[\ZZ]{L}{\QQ/\ZZ}}
  \end{equation*}

  is acyclic. The cycles of $\Hom[\ZZ]{L}{\QQ/\ZZ}$ have the form
  \begin{equation*}
    \Cy[i]{\Hom[\ZZ]{L}{\QQ/\ZZ}} \deq
    \Hom[\ZZ]{\Bo[-i-1]{L}}{\QQ/\ZZ}\:.
  \end{equation*}
  Now apply \prpref{ptwize-Z}.
\end{prf*}


\def\soft#1{\leavevmode\setbox0=\hbox{h}\dimen7=\ht0\advance \dimen7
  by-1ex\relax\if t#1\relax\rlap{\raise.6\dimen7
  \hbox{\kern.3ex\char'47}}#1\relax\else\if T#1\relax
  \rlap{\raise.5\dimen7\hbox{\kern1.3ex\char'47}}#1\relax \else\if
  d#1\relax\rlap{\raise.5\dimen7\hbox{\kern.9ex \char'47}}#1\relax\else\if
  D#1\relax\rlap{\raise.5\dimen7 \hbox{\kern1.4ex\char'47}}#1\relax\else\if
  l#1\relax \rlap{\raise.5\dimen7\hbox{\kern.4ex\char'47}}#1\relax \else\if
  L#1\relax\rlap{\raise.5\dimen7\hbox{\kern.7ex
  \char'47}}#1\relax\else\message{accent \string\soft \space #1 not
  defined!}#1\relax\fi\fi\fi\fi\fi\fi}
  \providecommand{\MR}[1]{\mbox{\href{http://www.ams.org/mathscinet-getitem?mr=#1}{#1}}}
  \renewcommand{\MR}[1]{\mbox{\href{http://www.ams.org/mathscinet-getitem?mr=#1}{#1}}}
  \providecommand{\arxiv}[2][AC]{\mbox{\href{http://arxiv.org/abs/#2}{\sf
  arXiv:#2 [math.#1]}}} \def\cprime{$'$}
\providecommand{\bysame}{\leavevmode\hbox to3em{\hrulefill}\thinspace}
\providecommand{\MR}{\relax\ifhmode\unskip\space\fi MR }
\providecommand{\MRhref}[2]{%
  \href{http://www.ams.org/mathscinet-getitem?mr=#1}{#2}
}
\providecommand{\href}[2]{#2}

\end{document}